\documentclass[11pt]{amsart}
\usepackage[margin=2.6cm]{geometry}
\usepackage{amssymb,amsmath,enumitem}
\usepackage[hidelinks]{hyperref}
\hypersetup{
  pdftitle={Alpha-points of a putative binary q-analog of the Fano plane},
  pdfauthor={Sascha Kurz},
  pdfsubject={A point P of a putative q-analog D of the Fano plane is an alpha-point if the derived design of D in P is a Desarguesian line spread. Kiermaier showed that for q prime or even every hyperplane contains a point which is not an alpha-point. For q=2 we study the 45 blocks of D inside a hyperplane H. Call a 5-subspace of H a beta-flat if all its blocks pass through one point, its focal point. We show that at least 41 points of H are not focal points of beta-flats. Consequently D has at least 85 non-alpha-points, at most 42 alpha-points, and at most 22 alpha-points in every hyperplane.},
  pdfkeywords={q-analog of the Fano plane, subspace design, alpha-point, hyperplane section, symplectic polarity}
}
\makeatletter
\@ifundefined{KV@Hyp@pdfcopyright}{}{%
  \hypersetup{pdfcopyright={Creative Commons Lizenz CC BY 4.0},pdflicenseurl={https://creativecommons.org}}}
\makeatother
\newtheorem{theorem}{Theorem}[section]
\newtheorem{lemma}[theorem]{Lemma}
\newtheorem{proposition}[theorem]{Proposition}
\newtheorem{corollary}[theorem]{Corollary}
\newtheorem{fact}[theorem]{Computational Fact}
\theoremstyle{definition}
\newtheorem{definition}[theorem]{Definition}
\newtheorem{remark}[theorem]{Remark}
\newtheorem{observation}[theorem]{Observation}
\newtheorem{problem}[theorem]{Problem}

\newcommand{\F}{\mathbb{F}}
\newcommand{\gauss}[3]{\genfrac{[}{]}{0pt}{}{#1}{#2}_{#3}}
\newcommand{\D}{\mathcal{D}}
\newcommand{\K}{\mathcal{K}}
\newcommand{\X}{\mathcal{X}}
\newcommand{\Y}{\mathcal{Y}}
\newcommand{\N}{\mathcal{N}}
\newcommand{\A}{\mathcal{A}}
\newcommand{\PG}{\operatorname{PG}}
\newcommand{\AG}{\operatorname{AG}}
\newcommand{\GL}{\operatorname{GL}}

\newcommand{\Der}{\operatorname{Der}}
\newcommand{\nnti}{n_{\mathrm{nTI}}}

\title[Alpha-points of a putative binary $q$-analog of the Fano plane]{$\text{\normalfont\large\bfseries Alpha-points of a putative binary $q$-analog of the Fano plane}$}
\author{$\text{\normalfont Sascha Kurz}$\\[3pt]
$\text{\normalfont\small Mathematisches Institut, Universit\"at Bayreuth, Universit\"atsstra\ss e 30, 95447 Bayreuth, Germany}$\\
$\text{\normalfont\small\texttt{sascha.kurz@uni-bayreuth.de}}$}
\date{}

\begin{document}
\begin{abstract}
A point $P$ of a putative $q$-analog $\D$ of the Fano plane is an $\alpha$-point if the derived design of $\D$ in $P$ is a Desarguesian line spread. Kiermaier showed that for $q$ prime or even every hyperplane contains a point which is not an $\alpha$-point. For $q=2$ we study the $45$ blocks of $\D$ inside a hyperplane $H$. Call a $5$-subspace of $H$ a $\beta$-flat if all its blocks pass through one point, its focal point. We show that at least $41$ points of $H$ are not focal points of $\beta$-flats. Consequently $\D$ has at least $85$ non-$\alpha$-points, at most $42$ $\alpha$-points, and at most $22$ $\alpha$-points in every hyperplane.

\medskip\noindent\textbf{Keywords:} $q$-analog of the Fano plane, subspace design, $\alpha$-point, hyperplane section, symplectic polarity.

\noindent\textbf{MSC classification (2020):} Primary 51E20, secondary 51E23, 05B25, 94B25.
\end{abstract}
\maketitle

\section{Introduction}

Let $V=\F_q^7$. A \emph{$q$-analog of the Fano plane} is a set $\D$ of $3$-subspaces (\emph{blocks}) of $V$ such that every $2$-subspace of $V$ lies in exactly one block, i.e.\ a $2$-$(7,3,1)_q$ subspace design. Its existence is open for every prime power $q$, see \cite{Kie22} and the references therein. The $q$-analog of Keevash's theorem \cite{Kee14,KSS25} gives $q$-Steiner systems $S_q(t,k,n)$ only for sufficiently large $n$, so it does not apply to $n=7$. The first nontrivial $q$-Steiner system, $S_2(2,3,13)$, was constructed in \cite{BEOVW16}. The question goes back to Thomas, who also excluded designs invariant under a Singer cycle of $\GL(7,2)$ \cite{Tho87}. For $q=2$ such a design would have $381$ blocks and its automorphism group would have order at most $2$ \cite{BKN16,KKW18}. The largest set of planes of $\PG(6,2)$ pairwise meeting in at most a point known to us has $333$ elements \cite{HKKW19}.

Equivalently, a $q$-analog of the Fano plane is a constant dimension code in the set of $3$-subspaces of $\F_q^7$ with minimum subspace distance $4$ and the largest possible size $q^8+q^6+q^5+q^4+q^3+q^2+1$. We write $A_q(7,4;3)$ for the maximum size of a set of planes of $\PG(6,q)$ pairwise meeting in at most a point. Finding large such sets is a set packing problem with one binary variable for each plane and one constraint for each line \cite{KK08}. For $q=2$ this standard formulation, with $11\,811$ binary variables, is the instance \texttt{cdc7-4-3-2} of MIPLIB~2017 \cite{GHG+21}, which is still open. By \cite{NS16} and \cite[Section~12]{Kur25}, $A_2(7,4;3)$ is neither $380$ nor $379$. So either a $q$-analog of the Fano plane exists or $A_2(7,4;3)\le378$.\footnote{The dual covering problem is open as well. Let $C_2(7,3,2)$ denote the minimum number of planes of $\PG(6,2)$ such that every line lies in at least one of them. The covering bound gives $C_2(7,3,2)\ge 381$, with equality if and only if a $q$-analog of the Fano plane exists. The best known upper bound is $C_2(7,3,2)\le 396$ \cite[Thm.~14]{Etz14}.}

For $q=3$ a $q$-analog of the Fano plane would have $7651$ blocks, and the best known lower bound is $A_3(7,4;3)\ge6978$. It comes from a code of size $6977$ \cite{HK16} to which one block can be added trivially (see \cite{HKKW16}, \cite[footnote~1]{Kie22}).

For a point $P$, the derived design $\Der_P(\D)=\{B/P: P\le B\in\D\}$ is a line spread of $V/P\cong\PG(5,q)$ \cite[Sect.~2.2]{Kie22}. Following \cite{HS16,Kie22}, $P$ is an \emph{$\alpha$-point} if this spread is geometric, equivalently Desarguesian \cite{BC74}. Thomas \cite{Tho96} proved that not all points are $\alpha$-points. Heden and Sissokho \cite{HS16} showed for $q=2$ that the non-$\alpha$-points meet every hyperplane. Kiermaier \cite{Kie22} introduced the tools on which our work rests: rich and poor solids, and $\beta$-flats, i.e.\ $5$-subspaces all of whose blocks pass through one point, their focal point. He proved for every $q$ that a hyperplane consisting of $\alpha$-points yields a partition of the generalized quadrangle $W(q)$ into spreads, which cannot exist for $q$ prime or even. This extends the result of Heden and Sissokho to all primes and all even $q$.

We ask instead how many points of a hyperplane $H$ must be non-$\alpha$-points, using only the $45$ blocks of $\D$ contained in $H$ (for $q=2$). Our starting point is Kiermaier's observation that every $\alpha$-point of $H$ is the focal point of a $\beta$-flat contained in $H$. Let $\X_H$ be the set of points of $H$ which are not focal points of $\beta$-flats in $H$, and let $\N$ and $\A$ denote the sets of non-$\alpha$-points and of $\alpha$-points of $\D$. Then $\X_H\subseteq\N\cap H$, and lower bounds on $|\X_H|$ transfer to global bounds via the Griesmer bound. This differs from the point-based viewpoint of Etzion and Hooker \cite{EH18}, who consider the derived and the residual design of a point. We use instead the blocks contained in a hyperplane, which form the residual design in the sense of Kiermaier and Laue \cite{KL15}.

\smallskip\noindent\textbf{General $q$.} Every hyperplane of a $q$-analog of the Fano plane, $q$ prime or even, contains a $5$-subspace which is not a $\beta$-flat (Theorem~\ref{thm:A}). This sharpens Kiermaier's result with the same method. It is our only result for general $q$ (Section~\ref{sec:general}).

\smallskip\noindent\textbf{The binary case $q=2$.}
\begin{enumerate}[label=(\roman*),leftmargin=*]
\item Every hyperplane contains at least $41$ points which are not focal points of $\beta$-flats in $H$, i.e.\ $|\X_H|\ge41$ (Theorem~\ref{thm:41}). This gives $|\N|\ge85$ and $|\A|\le42$, at most $22$ $\alpha$-points in every hyperplane, and at least $1736$ of the $2667$ five-subspaces are not $\beta$-flats (Corollary~\ref{cor:41}). The proof splits according to whether the set $\A'$ of focal points of the $\beta$-flats in $H$ spans $H$. If it does, the $\beta$-flats induce a symplectic polarity on $H$, and a parity argument gives even $|\X_H|\ge44$ (Lemma~\ref{lem:polspan}, Theorem~\ref{thm:pol44}). If it does not, a counting argument on the blocks outside the $5$-subspace spanned by $\A'$ excludes $|\X_H|\le40$ (Section~\ref{sec:41}). The proof is by hand, and computer enumerations serve as cross-checks (Appendix~\ref{app:repro}).
\item All bounds on $|\X_H|$ hold for abstract \emph{hyperplane sections}, i.e.\ sets of $45$ planes of $\PG(5,2)$ pairwise meeting in at most a point with exactly $5$ through every point (Section~\ref{sec:dh}). There we also study the lifting problem and show, by computer \textup{(}Fact~\ref{fact:w}\textup{)}, that the linear relations between the numbers of $5$-subspace patterns do not exclude a Fano design.
\end{enumerate}

The remaining part of the paper is structured as follows.
Section~\ref{sec:general} treats general $q$, from Section~\ref{sec:prelim} on we assume $q=2$, Sections~\ref{sec:prelim}--\ref{sec:41} prove Theorem~\ref{thm:41}, and Section~\ref{sec:dh} studies lifting. We draw a brief conclusion and present some open problems in Section~\ref{sec:open}.  Appendix~\ref{app:repro} describes the programs and cross-checks, and Appendix~\ref{app:section} lists an explicit hyperplane section.

\section{\texorpdfstring{Results for general $q$}{Results for general q}}\label{sec:general}
In this section $q$ is an arbitrary prime power, unless stated otherwise. From Section~\ref{sec:prelim} on we restrict to $q=2$. Let $V=\F_q^7$ and let $\D$ be a $2$-$(7,3,1)_q$ design, i.e.\ a $q$-analog of the Fano plane.

\subsection{\texorpdfstring{Designs, spreads and $\alpha$-points}{Designs, spreads and alpha-points}}
We use algebraic dimensions: points, lines, planes and solids are subspaces of dimension $1,2,3,4$, and a hyperplane of $V=\F_q^7$ has dimension $6$. $\gauss nkq$ is the Gaussian binomial coefficient. A \emph{line spread} of a vector space is a set of lines partitioning its points. For a $2$-$(7,3,1)_q$ design $\D$, a point $P$, a hyperplane $H$ and a $5$-subspace $F$ one has \cite[Sect.~2.2]{Kie22}
\[
|\D|=q^8+q^6+q^5+q^4+q^3+q^2+1,\quad \#\{B\ni P\}=q^4+q^2+1,\quad \#\{B: P\le B\le H\}=q^2+1,
\]
\[
\#\{B: B\le H\}=q^5+q^3+q^2+1,\qquad \#\{B: B\le F\}=q^2+1.
\]
Two distinct blocks meet in a point or trivially. So two distinct blocks in a common $5$-subspace meet in exactly one point, and a solid contains at most one block. A solid is \emph{rich} if it contains a block and \emph{poor} otherwise \cite{Kie22}.

A line spread $\mathcal S$ of $\F_q^6$ is \emph{geometric} if every solid spanned by two of its lines contains $q^2+1$ of its lines. Geometric spreads are exactly the Desarguesian ones \cite{BC74}, \cite[Sect.~2.3]{Kie22}. A point $P$ is an \emph{$\alpha$-point} of $\D$ if $\Der_P(\D)$ is geometric.

\begin{definition}[\cite{Kie22}]
A $5$-subspace $F$ is a \emph{$\beta$-flat with focal point $P$} if all $q^2+1$ blocks contained in $F$ pass through $P$.
\end{definition}

Heuristically, a $\beta$-flat with focal point $P$ plays the role of $P^\perp$ for a symplectic polarity. Theorem~\ref{thm:A} shows that this is literally true if all $5$-subspaces of $H$ are $\beta$-flats.

\begin{lemma}[Kiermaier {\cite[Lemmas~3.2--3.5, 3.8]{Kie22}}]\label{lem:K}
Let $\D$ be a $2$-$(7,3,1)_q$ design and $H$ a hyperplane.
\begin{enumerate}[label=(\alph*)]
\item The focal point of a $\beta$-flat is unique, and a point $P$ of $H$ is the focal point of at most one $\beta$-flat contained in $H$.
\item If $F$ is a $\beta$-flat with focal point $P$, the blocks in $F$ partition the points of $F$ other than $P$. The poor solids of $F$ are exactly the solids of $F$ not through $P$, and each poor solid $S$ of $F$ is partitioned into lines by $\{B\cap S: B\in\D,\ B\le F\}$.
\item If $P$ is an $\alpha$-point and $B\cap B'=P$ for blocks $B,B'$, then $B+B'$ is a $\beta$-flat with focal point $P$. In particular, every $\alpha$-point $P\in H$ is the focal point of a $\beta$-flat contained in $H$.
\end{enumerate}
\end{lemma}

\begin{theorem}[Kiermaier {\cite[Thms.~1, 2]{Kie22}}, for $q=2$ Heden--Sissokho \cite{HS16}]\label{thm:K}
If some hyperplane consists of $\alpha$-points, then the line set of $W(q)$ can be partitioned into spreads. If $q$ is prime or even, this is impossible, so every hyperplane contains a non-$\alpha$-point.
\end{theorem}

Here $W(q)$ is the generalized quadrangle formed by the points of $\PG(3,q)$ and the totally isotropic lines of a nondegenerate alternating form, and a spread of $W(q)$ is a set of its lines partitioning the points \cite[Sect.~1.1, 1.8 and 3.1.1]{PT09}.

By Lemma~\ref{lem:K}(a),(c), the $\gauss61q$ points of a hyperplane $H$ consisting of $\alpha$-points are focal points of pairwise distinct $\beta$-flats contained in $H$, and $H$ has only $\gauss65q=\gauss61q$ five-subspaces. So all of them are $\beta$-flats, and the following theorem contains the last statement of Theorem~\ref{thm:K}.

\begin{theorem}\label{thm:A}
Let $\D$ be a $q$-analog of the Fano plane \textup{(}$q$ arbitrary\textup{)} and $H$ a hyperplane. If every $5$-subspace of $H$ is a $\beta$-flat, then the line set of $W(q)$ can be partitioned into spreads. Consequently, if $q$ is prime or even, every hyperplane contains a $5$-subspace which is not a $\beta$-flat.
\end{theorem}
\begin{proof}
Let $G(Q)$ denote the $\beta$-flat with focal point $Q$. We show that $Q\mapsto G(Q)$ is a symplectic polarity of $H$, and that the $\beta$-flats through a nondegenerate solid yield $q+1$ pairwise disjoint spreads of $W(q)$.

By Lemma~\ref{lem:K}(a) and $\gauss65q=\gauss61q$, the map sending a $5$-subspace of $H$ to its focal point is a bijection onto the points of $H$. By Lemma~\ref{lem:K}(b), the $q^2+1$ blocks in $G(Q)$ are the blocks contained in $H$ through $Q$, and they cover $G(Q)$. Hence for points $Q\ne R$ of $H$ we have $R\in G(Q)$ iff some block contained in $H$ contains $Q$ and $R$ iff $Q\in G(R)$. Also $Q\in G(Q)$.

Let $Q_1,Q_2,Q_3$ be collinear and $R\in G(Q_1)\cap G(Q_2)$. Then $Q_1,Q_2\in G(R)$, so $Q_3\in G(R)$ and $R\in G(Q_3)$. Thus $G(Q_3)$ contains the solid $G(Q_1)\cap G(Q_2)$, and $Q\mapsto G(Q)$ maps the points of a line bijectively onto the $5$-subspaces through a solid. So $Q\mapsto G(Q)$ is a collineation from $\PG(H)$ onto its dual. By the fundamental theorem of projective geometry \cite{Tay92} there are an automorphism $\sigma$ of $\F_q$ and a $\sigma$-semilinear bijection $\varphi:H\to H^*$ with $G(\langle u\rangle)=\ker\varphi(u)$. Put $f(u,v)=\varphi(u)(v)$. Since $Q\in G(Q)$, $f(u,u)=0$ for all $u$, and expanding $f(u+v,u+v)$ gives $f(u,v)=-f(v,u)$. Then $\lambda^\sigma f(u,v)=f(\lambda u,v)=-f(v,\lambda u)=\lambda f(u,v)$, so $\sigma=1$ as $f\ne0$. Hence $f$ is a nondegenerate alternating form with $G(Q)=Q^\perp$ for all points $Q$.

Every block $B\le H$ lies in $G(P)=P^\perp$ for each point $P\in B$, so $B$ is totally isotropic. Choose points $Q,Q'$ with $f(Q,Q')\ne0$, let $\ell=Q+Q'$ and $S=\ell^\perp$. Then $S\cap\ell=0$ and $f$ is nondegenerate on $S$. As $Q\notin S\le G(Q)$, $S$ is a poor solid of the $\beta$-flat $G(Q)$. The $q+1$ five-subspaces of $H$ through $S$ are $G(R)=R^\perp$, $R\in\ell$, and by Lemma~\ref{lem:K}(b), $T_R=\{B\cap S:B\le G(R)\}$ is a set of $q^2+1$ lines partitioning $S$. These lines are totally isotropic, so $T_R$ is a spread of the copy of $W(q)$ defined by $f|_S$. A line in $T_R\cap T_{R'}$ would lie in a block $B$ with $B+S=G(R)=G(R')$, so the $T_R$ are pairwise disjoint. Since $W(q)$ has $(q+1)(q^2+1)$ lines, they partition its line set. For $q$ prime or even this is impossible by \cite[proof of Thm.~2]{Kie22}, and for $q=2$ see also Lemma~\ref{lem:W2}.
\end{proof}

Let $\X_H$ be the set of points of a hyperplane $H$ which are not focal points of $\beta$-flats contained in $H$, and let $\N$ be the set of non-$\alpha$-points of $\D$. By Lemma~\ref{lem:K}(c), every $\alpha$-point of $H$ is the focal point of a $\beta$-flat contained in $H$, so $\X_H\subseteq\N\cap H$.

\section{\texorpdfstring{The binary case: Fano designs and hyperplane sections}{The binary case: Fano designs and hyperplane sections}}\label{sec:prelim}
From now on $q=2$, and the letter $q$ is no longer used as a variable. A \emph{Fano design} is a $2$-$(7,3,1)_2$ design, i.e.\ a binary $q$-analog of the Fano plane. Specializing the numbers of Section~\ref{sec:general}: a Fano design $\D$ has $|\D|=381$ blocks, every point lies on $21$ blocks, every hyperplane $H$ contains $45$ blocks, every point of $H$ lies on $5$ blocks contained in $H$, and every $5$-subspace contains $5$ blocks.

All our bounds use only the blocks contained in one hyperplane, and we formulate them for the following abstract objects.

\begin{definition}
Let $H$ be a $6$-dimensional vector space over $\F_2$. A \emph{hyperplane section} on $H$ is a set $\K$ of planes of $H$ (\emph{blocks}) such that no line of $H$ lies in two blocks and every point of $H$ lies on exactly $5$ blocks. A line of $H$ is a \emph{block line} if it lies in a block and \emph{free} otherwise. A $5$-subspace $F\le H$ is a \emph{$\beta$-flat with focal point $P$} if all blocks contained in $F$ pass through $P$. A solid is \emph{rich} if it contains a block and \emph{poor} otherwise.
\end{definition}

In the language of subspace designs, a hyperplane section is a $1$-$(6,3,5)_2$ design which is also a constant dimension code of minimum subspace distance $4$ \textup{(}a partial $2$-$(6,3,1)_2$ design\textup{)}. An explicit hyperplane section with $15$ $\beta$-flats is given in Appendix~\ref{app:section}.

If $\D$ is a Fano design and $H$ a hyperplane of $V$, then $\D|_H=\{B\in\D: B\le H\}$ is a hyperplane section on $H$, with the same $\beta$-flats.

Parts (b) and (c) of the following lemma are Lemma~\ref{lem:K}(a),(b) for $q=2$. Since we need them for abstract hyperplane sections, which do not come from a design, we include the short proofs.

\begin{lemma}[cf.\ {\cite[Lemmas~3.2--3.4]{Kie22}}]\label{lem:dh}
Let $\K$ be a hyperplane section on $H$.
\begin{enumerate}[label=(\alph*)]
\item $|\K|=45$, with $315$ block lines and $336$ free lines. Every $5$-subspace $F\le H$ contains exactly $5$ blocks, which pairwise meet in exactly one point, and every other block meets $F$ in a line. A solid contains at most one block.
\item The focal point of a $\beta$-flat is unique, and every point is the focal point of at most one $\beta$-flat.
\item If $F$ is a $\beta$-flat with focal point $P$, then the $5$ blocks of $F$ are the $5$ blocks through $P$, they partition $F\setminus\{P\}$, the poor solids of $F$ are exactly the $16$ solids of $F$ not through $P$, and every poor solid $S$ of $F$ is partitioned into the lines $B\cap S$, $B\le F$.
\end{enumerate}
\end{lemma}
\begin{proof}
(a) Counting incidences of points and blocks gives $7|\K|=5\cdot63$, so $|\K|=45$, and $7\cdot45=315$ of the $651$ lines are block lines. Two blocks meet in at most a point, so a solid contains at most one block, and a block not in $F$ meets $F$ in a line. If $F$ contains $b$ blocks, counting incidences with the $31$ points of $F$ gives $5\cdot31=7b+3(45-b)$, so $b=5$. Two blocks of $F$ meet in a point by the dimension formula.

(b) Two focal points $P\ne Q$ would put the line $P+Q$ into $5$ blocks. A $\beta$-flat with focal point $P$ contains $5$ blocks through $P$, i.e.\ all blocks through $P$. Two such $\beta$-flats would meet in a solid containing $5$ blocks.

(c) The $5$ blocks through $P$ pairwise meet only in $P$ and cover $1+5\cdot6=31$ points. Each block lies in $3$ solids of $F$, so the $15$ rich solids of $F$ all pass through $P$, and $F$ has exactly $15$ solids through $P$ and $16$ solids not through $P$. If $P\notin S$, each $B\cap S$ is a line, and these lines partition $S$.
\end{proof}

\begin{definition}
For a hyperplane section on $H$ let $\Y=\Y_H$ be the set of $5$-subspaces of $H$ which are not $\beta$-flats, $\A'=\A'_H$ the set of focal points of $\beta$-flats, and $\X=\X_H=H\setminus\A'$. For $Q\in\A'$ let $G(Q)$ be the $\beta$-flat with focal point $Q$. For a $5$-subspace $F$ and $P\in F$ let $j_F(P)$ be the number of blocks of $F$ through $P$.
\end{definition}

By Lemma~\ref{lem:dh}(b), $F\mapsto$ focal point is a bijection from the $\beta$-flats onto $\A'$, so $|\X_H|=|\Y_H|$ (there are $63$ points and $63$ five-subspaces in $H$). For a hyperplane section $\D|_H$ of a Fano design $\D$, the set $\X$ is the set $\X_H$ of Section~\ref{sec:general}.

All statements of Sections~\ref{sec:tools}--\ref{sec:41} about $\X_H$ are proved for arbitrary hyperplane sections. Only the passage to $\N$ and $\A$ uses the design (Lemma~\ref{lem:K}(c)).

\begin{remark}
Lemma~\ref{lem:dh} and the counting statements of Section~\ref{sec:dh} (Propositions~\ref{prop:lift} and \ref{prop:pencil}) should have analogues for general $q$, but we restrict to $q=2$. The construction of the form in Lemma~\ref{lem:polspan} and the arguments of Lemmas~\ref{lem:parity} and \ref{lem:W2} use $q=2$.
\end{remark}

\section{\texorpdfstring{Block patterns, poor solids and spreads of $W(2)$}{Block patterns, poor solids and spreads of W(2)}}\label{sec:tools}

For a $5$-subspace $F$ of a hyperplane section, let $m_j(F)$ be the number of points $P\in F$ with $j_F(P)=j$. The vector $\mathbf m(F)=(m_0,\dots,m_5)$ is the \emph{pattern} of $F$.

\begin{proposition}\label{prop:patterns}
Let $B_1,\dots,B_5$ be planes of $\F_2^5$ pairwise meeting in exactly one point. Then $(m_0,\dots,m_5)$ is one of
\begin{gather*}
t_1=(0,30,0,0,0,1),\quad t_2=(3,23,4,0,1,0),\quad t_3=(4,21,4,2,0,0),\\ t_4=(5,18,7,1,0,0),\quad t_5=(6,15,10,0,0,0).
\end{gather*}
The number $\eta(F)=m_0+m_2+m_3+m_4$ of points of $F$ on $0$, $2$, $3$ or $4$ of the planes is $0,8,10,13,16$, respectively, and $t_1$ occurs iff the five planes have a common point.
\end{proposition}
\begin{proof}
Double counting gives
\[
\sum_jm_j=31,\qquad \sum_jj\,m_j=35,\qquad \sum_j\binom j2m_j=\binom52=10,
\]
so the pattern is determined by $(m_3,m_4,m_5)$.

For a point $X$ let $J(X)=\{i:X\in B_i\}$. If $X\ne Y$ and $|J(X)\cap J(Y)|\ge2$, two planes would share the line $X+Y$, so the sets $J(X)$ with $|J(X)|\ge2$ are the blocks of a linear space on $\{1,\dots,5\}$ (every pair of indices lies in exactly one of them, since two planes meet in exactly one point), and $m_j$ is its number of blocks of size $j$ for $j\ge2$. A block of size $5$ or $4$ is the only block of size $\ge3$ (a $4$-subset meets every other $3$-subset of a $5$-set in at least two elements). Two $3$-blocks meet in one element and cover all five elements, so there is no third. Hence $(m_3,m_4,m_5)\in\{(0,0,1),(0,1,0),(2,0,0),(1,0,0),(0,0,0)\}$, which gives $t_1,\dots,t_5$.
\end{proof}

A $5$-subspace $F$ is a $\beta$-flat iff $\mathbf m(F)=t_1$.

\begin{lemma}\label{lem:focal}
Let $Q$ be the focal point of a $\beta$-flat $G\le H$ and $F\ne G$ a $5$-subspace of $H$ with $Q\in F$. Then $j_F(Q)=1$. So if $F\in\Y$, every point $P\in F$ with $j_F(P)\ne1$ lies in $\X$.
\end{lemma}
\begin{proof}
The blocks of $F$ through $Q$ are blocks of $G$ (Lemma~\ref{lem:dh}(c)) and lie in the solid $F\cap G$, so $j_F(Q)\le1$. Each of the $5$ blocks through $Q$ lies in $7$ five-subspaces of $H$, so $\sum_{F\ni Q}j_F(Q)=35$, where $F$ runs over the $31$ five-subspaces of $H$ through $Q$. Here $G$ contributes $5$ and each of the other $30$ at most $1$, hence exactly $1$. If $F\in\Y$ and $P\in F\cap\A'$, then $G(P)\ne F$ and $j_F(P)=1$.
\end{proof}

From here on in this section $\K$ is a hyperplane section on $H$, with $\X$ and $\Y$ as in Section~\ref{sec:prelim}.

\begin{lemma}\label{lem:forced}
Let $S\le H$ be a poor solid, $F_1,F_2,F_3$ the $5$-subspaces with $S<F_i<H$, and $T_i=\{B\cap S: B\in\K,\ B\le F_i\}$.
\begin{enumerate}[label=(\roman*)]
\item Each $T_i$ is a set of $5$ lines of $S$, and $T_1,T_2,T_3$ are pairwise disjoint.
\item If $F_i$ is a $\beta$-flat, then $T_i$ is a line spread of $S$.
\end{enumerate}
\end{lemma}
\begin{proof}
(i) A block $B\le F_i$ is not contained in $S$, and $\dim(B\cap S)\ge3+4-5$, so $B\cap S$ is a line and $B+S=F_i$. Distinct blocks give distinct lines, and a line in $T_i\cap T_j$ would lie in a block $B$ with $B+S=F_i=F_j$.

(ii) Lemma~\ref{lem:dh}(c).
\end{proof}

\begin{lemma}\label{lem:spreadplane}
Every plane of $\PG(3,2)$ contains exactly one line of a line spread $T$ of $\PG(3,2)$.
\end{lemma}
\begin{proof}
$\PG(3,2)$ has $15$ planes, with $3$ planes through each of its $35$ lines. Two lines of a plane meet, so a plane contains at most one line of $T$. The $5$ lines of $T$ lie in $5\cdot3=15$ planes, so every plane contains exactly one.
\end{proof}

\begin{lemma}\label{lem:W2}
$W(2)$ has exactly $6$ spreads, and any two of them share exactly one line.
\end{lemma}
\begin{proof}
This is the classical description of $W(2)$ by the duads and synthemes of a $6$-set, whose six spreads are the synthematic totals \textup{(}see \cite{PT09}\textup{)}. We give a direct proof. Let $\langle\cdot,\cdot\rangle$ be a nondegenerate alternating form on $V=\F_2^4$. The totally isotropic lines through a point $x$ are the three lines through $x$ in the plane $x^\perp$, so $W(2)$ has $15$ lines. Fix a non-TI line $h$. Then $V=h\oplus h^\perp$, and the form is nondegenerate on $h$ and on $h^\perp$, so $\langle x,x'\rangle=1$ for distinct points $x,x'$ of $h$, and likewise in $h^\perp$. For a point $x$ of $h$ we have $x^\perp=x+h^\perp$, so the TI lines through $x$ are the \emph{grid lines} $g(x,y)=\{x,y,x+y\}$ with $y$ a point of $h^\perp$. For a bijection $\sigma$ from the points of $h$ to the points of $h^\perp$ let $L_\sigma=\{x+\sigma(x):x\in h\}$. This is a line, since the three points sum to $0$, and it is TI, since $\langle x+y,x'+y'\rangle=\langle x,x'\rangle+\langle y,y'\rangle=1+1=0$ for $x\ne x'$ and $y\ne y'$. The $9$ grid lines and the $6$ lines $L_\sigma$ are distinct, so these are all $15$ lines of $W(2)$.

For bijections $\pi,\sigma$ from the points of $h$ to the points of $h^\perp$, $\pi^{-1}\sigma$ is a permutation of the three points of $h$. Let $C$ be the alternating group on these three points. Let $S$ be a spread. The lines of $S$ through the points $x$ of $h$ are grid lines $g(x,\pi(x))$, where $\pi$ is a bijection since the lines of $S$ are disjoint. These three lines cover $h$, $h^\perp$ and the points $x+\pi(x)$. The other two lines of $S$ miss $h\cup h^\perp$, so they are of the form $L_\sigma,L_\tau$, and for every $x$ the points $x+\pi(x)$, $x+\sigma(x)$, $x+\tau(x)$ are distinct. Thus $\pi^{-1}\sigma$, $\pi^{-1}\tau$ and $\sigma^{-1}\tau$ are fixed-point-free permutations of the points of $h$, i.e.\ $\{\sigma,\tau\}=\pi C\setminus\{\pi\}$ (the coset taken with respect to composition on the right, i.e.\ $\sigma=\pi\circ c$ with $c\in C$). Conversely, $S_\pi=\{g(x,\pi(x)):x\in h\}\cup\{L_\sigma:\sigma\in\pi C\setminus\{\pi\}\}$ is a spread. So the spreads are the $6$ sets $S_\pi$.

To list them, label the points of $h$ as $x_1,x_2,x_3$ and those of $h^\perp$ as $y_1,y_2,y_3$, and write a bijection $\sigma$ as the word $\sigma(x_1)\sigma(x_2)\sigma(x_3)$ of indices, so $231$ means $x_1\mapsto y_2$, $x_2\mapsto y_3$, $x_3\mapsto y_1$. Then $\pi C$ consists of the cyclic rotations of the word of $\pi$, and the six spreads are
\begin{center}\small
\begin{tabular}{@{}lll@{}}
$\pi$ & grid lines $g(x_i,y_{\pi(i)})$ & lines $L_\sigma$\\\hline
$123$ & $g(x_1,y_1),\,g(x_2,y_2),\,g(x_3,y_3)$ & $L_{231},\,L_{312}$\\
$231$ & $g(x_1,y_2),\,g(x_2,y_3),\,g(x_3,y_1)$ & $L_{123},\,L_{312}$\\
$312$ & $g(x_1,y_3),\,g(x_2,y_1),\,g(x_3,y_2)$ & $L_{123},\,L_{231}$\\
$132$ & $g(x_1,y_1),\,g(x_2,y_3),\,g(x_3,y_2)$ & $L_{321},\,L_{213}$\\
$321$ & $g(x_1,y_3),\,g(x_2,y_2),\,g(x_3,y_1)$ & $L_{132},\,L_{213}$\\
$213$ & $g(x_1,y_2),\,g(x_2,y_1),\,g(x_3,y_3)$ & $L_{132},\,L_{321}$
\end{tabular}
\end{center}

Let $\pi\ne\rho$. If $\pi^{-1}\rho\in C$, it has no fixed point, so $S_\pi$ and $S_\rho$ share no grid line. Since $\pi C=\rho C$, they share exactly the line $L_\tau$ with $\tau$ the third element of $\pi C$. Otherwise $\pi^{-1}\rho$ is a transposition. It has exactly one fixed point, which gives exactly one common grid line, and $\pi C\ne\rho C$ gives no common line $L_\sigma$.
\end{proof}

\section{\texorpdfstring{The polarity induced by the $\beta$-flats}{The polarity induced by the beta-flats}}\label{sec:polar}

Let $\K$ be a hyperplane section on $H$. For $Q\in\A'$ let $a_Q\in H^*$ be the linear form with $\ker a_Q=G(Q)$. For a point $P$ let $N(P)$ be the union of the $5$ blocks through $P$ ($31$ points).

\begin{lemma}\label{lem:GQ}
For $Q\in\A'$, $G(Q)=N(Q)$. So for $Q\in\A'$ and a point $R\ne Q$, $R\in G(Q)$ iff $Q$ and $R$ lie on a common block iff $Q\in N(R)$, and for $R\in\A'$ this reads $a_Q(R)=a_R(Q)$.
\end{lemma}
\begin{proof}
$G(Q)$ consists of the $5$ blocks through $Q$, which pairwise meet only in $Q$ and cover $31$ points (Lemma~\ref{lem:dh}(c)).
\end{proof}

\begin{lemma}\label{lem:blocksid}
A block $B$ lies in exactly $7-|B\cap\X|$ $\beta$-flats \textup{(}namely $G(Q)$, $Q\in B\setminus\X$\textup{)} and in exactly $|B\cap\X|$ members of $\Y$. Moreover $\sum_B|B\cap\X|=5|\X_H|$.
\end{lemma}
\begin{proof}
A $\beta$-flat $G(Q)\ge B$ has all its blocks through $Q$, so $Q\in B$. Conversely $B\le G(Q)$ for $Q\in B\setminus\X$ by Lemma~\ref{lem:GQ}, and $Q\mapsto G(Q)$ is injective. $B$ lies in $7$ five-subspaces of $H$. Each point lies on $5$ blocks, so $\sum_B|B\cap\X|=5|\X_H|$.
\end{proof}

\begin{lemma}\label{lem:polspan}
Suppose that $\A'$ spans $H$. Then there is an alternating bilinear form $\langle\cdot,\cdot\rangle$ on $H$ with $Q^\perp=G(Q)$ for every $Q\in\A'$. Its radical $R$ satisfies $R\cap\A'=\emptyset$ and $\dim R\in\{0,2\}$. If $\dim R=2$, then $|\A'|\le15$. In particular, if $\A'$ spans $H$ and $|\X_H|\le47$, the form is nondegenerate.
\end{lemma}
\begin{proof}
We construct the form from the linear forms $a_Q$ and then determine its radical. For $u\in H$ write $u=\sum_iQ_i$ with $Q_i\in\A'$ and put $M(u)=\sum_ia_{Q_i}\in H^*$. This is well defined: if $\sum_iQ_i=\sum_jQ'_j$, then for every $R\in\A'$, by Lemma~\ref{lem:GQ}, $\sum_ia_{Q_i}(R)=a_R(\sum_iQ_i)=a_R(\sum_jQ'_j)=\sum_ja_{Q'_j}(R)$, and linear forms which agree on the spanning set $\A'$ are equal. So $M:H\to H^*$ is linear.

Put $\langle u,v\rangle=M(u)(v)$. For $Q,R\in\A'$, $\langle Q,R\rangle=a_Q(R)=a_R(Q)=\langle R,Q\rangle$, so the form is symmetric. Since the form is symmetric and the characteristic is $2$, the cross terms cancel, so $\langle u,u\rangle=\sum_i\langle Q_i,Q_i\rangle=\sum_ia_{Q_i}(Q_i)=0$, because $Q_i\in G(Q_i)$. Hence the form is alternating. 

For $Q\in\A'$, $Q^\perp=\ker a_Q=G(Q)$, and $G(Q)\ne H$, so $Q\notin R$ and $R\cap\A'=\emptyset$. The rank of an alternating form is even, so $\dim R\in\{0,2,4,6\}$, and $\dim R=6$ is impossible since it would give $Q\in R$. For $Q\in\A'$ we have $G(Q)=Q^\perp\supseteq R$, and $Q\mapsto G(Q)$ is injective. A solid lies in $3$ and a line in $15$ five-subspaces of $H$, so $\dim R=4$ gives $|\A'|\le3$, which contradicts the spanning assumption, and $\dim R=2$ gives $|\A'|\le15$. If $|\X_H|\le47$, then $|\A'|=63-|\X_H|\ge16$ and $R=0$.
\end{proof}

We say that \emph{$\K$ carries the polarity} if there is a nondegenerate alternating form on $H$ with $G(Q)=Q^\perp$ for all $Q\in\A'$. A subspace is \emph{totally isotropic} (TI) if the form vanishes on it.

We use the following facts on a nondegenerate alternating form on $H\cong\F_2^6$, which are also confirmed by computer (Appendix~\ref{app:repro}). A plane $B$ is TI or has rank $2$. Indeed, $\dim B^\perp=3$ and the rank of the form on $B$ is even, so either $B\le B^\perp$ or $B\cap B^\perp$ is a single point $r_B$, the unique radical point of $B$. Let $\nnti$ be the number of non-TI blocks of $\K$.

\begin{lemma}\label{lem:polcons}
Suppose that $\K$ carries the polarity. Every block through a point $Q\in\A'$ lies in $Q^\perp$. For every non-TI block $B$, $(B\cup B^\perp)\setminus\{r_B\}\subseteq\X$, so $B\cap\A'\subseteq\{r_B\}$ and $|B\cap\X|\ge6$.
\end{lemma}
\begin{proof}
$B\le G(Q)=Q^\perp$ by Lemma~\ref{lem:dh}(c). Let $B$ be non-TI. If $u\in B\cap\A'$, then $B\le u^\perp$, i.e.\ $u\in B\cap B^\perp=\{r_B\}$. If $u\in B^\perp\cap\A'$, then $B\le u^\perp=G(u)$, so $B$ is a block of the $\beta$-flat $G(u)$ and passes through $u$. Again $u=r_B$.
\end{proof}

\begin{lemma}\label{lem:parity}
Suppose that $\K$ carries the polarity. For $Q\in\A'$ let $\nu(Q)$ be the number of non-TI blocks through $Q$. Then $\nu(Q)\in\{0,2,4\}$.
\end{lemma}
\begin{proof}
Let $\bar V=Q^\perp/Q$, a $4$-dimensional space. The radical of the restriction of the form to $Q^\perp$ is $Q^\perp\cap(Q^\perp)^\perp=Q^\perp\cap Q=Q$, so $\langle u+Q,v+Q\rangle=\langle u,v\rangle$ is a nondegenerate alternating form on $\bar V$. By Lemma~\ref{lem:dh}(c), the five blocks through $Q$ lie in $G(Q)=Q^\perp$ and partition $Q^\perp\setminus\{Q\}$. Each of them is a union of cosets of $Q$, so their images $B/Q$ are five lines partitioning the $15$ points of $\bar V$. As $Q$ is orthogonal to $B$, the block $B$ is TI iff $B/Q$ is TI. Hence $\nu(Q)$ is the number of non-TI lines of this partition.

Let $\Pi$ be the set of unordered pairs $\{x,y\}$ of points of $\bar V$ with $\langle x,y\rangle=1$. For every point $x$ there are $8$ points $y$ with $\langle x,y\rangle=1$, so $|\Pi|=15\cdot8/2=60$. Two distinct points of a non-TI line span it, so all $3$ pairs of points of a non-TI line lie in $\Pi$, and no pair of points of a TI line does. If $L\ne M$ are lines of the partition and $x\in L$, then $y\mapsto\langle x,y\rangle$ is a linear form on $M$, which vanishes on all of $M$ or takes the value $1$ on exactly two of its three points. So an even number of pairs in $\Pi$ have one point on $L$ and one on $M$. Therefore $60\equiv3\nu(Q)\pmod2$, and $\nu(Q)\le5$ is even.
\end{proof}

\begin{lemma}\label{lem:A0}
Suppose that $\K$ carries the polarity, and let $\A'_0=\{Q\in\A':\nu(Q)=0\}$. Then $\langle Q,Q'\rangle=0$ for all $Q,Q'\in\A'_0$. Hence $\langle\A'_0\rangle$ is totally isotropic, $\dim\langle\A'_0\rangle\le3$ and $|\A'_0|\le7$.
\end{lemma}
\begin{proof}
Let $Q\ne Q'$ in $\A'_0$ with $\langle Q,Q'\rangle=1$, let $\ell=Q+Q'$ and $S=\ell^\perp=G(Q)\cap G(Q')$. This is a solid, since $G(Q)\ne G(Q')$ (Lemma~\ref{lem:dh}(b)). As $\langle Q,Q'\rangle=1$, neither $Q$ nor $Q'$ lies in $S$. So $S$ is a solid of the $\beta$-flat $G(Q)$ not through its focal point, hence poor (Lemma~\ref{lem:dh}(c)). The line $\ell$ is not TI, so $S\cap S^\perp=\ell^\perp\cap\ell=0$, and the form restricted to $S$ is nondegenerate. The flats $G(Q)$ and $G(Q')$ are two of the three $5$-subspaces through $S$. By Lemma~\ref{lem:forced}(i),(ii), $T_Q=\{B\cap S:Q\in B\in\K\}$ and $T_{Q'}$ are disjoint line spreads of $S$. Since $\nu(Q)=\nu(Q')=0$, all blocks through $Q$ or $Q'$ are TI, so all lines of $T_Q$ and $T_{Q'}$ are totally isotropic. Hence $T_Q$ and $T_{Q'}$ are disjoint spreads of the generalized quadrangle $W(2)$ defined by the restricted form on $S$, which contradicts Lemma~\ref{lem:W2}. Finally, pairwise orthogonal points span a totally isotropic subspace $W$, and $W\le W^\perp$ gives $\dim W\le3$.
\end{proof}

\begin{theorem}\label{thm:pol44}
If a hyperplane section carries the polarity, then $|\X_H|\ge44$.
\end{theorem}
\begin{proof}
Idea: every point of $\A'\setminus\A'_0$ lies on at least two non-TI blocks (Lemma~\ref{lem:parity}), and every non-TI block uses six points of $\X$ (Lemma~\ref{lem:polcons}). The points of $\A'_0$ are handled by counting TI blocks.

Let $x=|\X_H|$, $s=|\A'|=63-x$, $a_0=|\A'_0|$ and let $n_{\mathrm{TI}}=45-\nnti$ be the number of TI blocks. By Lemma~\ref{lem:polcons}, a non-TI block $B$ satisfies $B\cap\A'\subseteq\{r_B\}$ and $|B\cap\X|\ge6$. Counting the pairs $(Q,B)$ with $Q\in\A'\cap B$ and $B$ non-TI, and using Lemma~\ref{lem:parity}, we get
\[
2(s-a_0)\le\sum_{Q\in\A'}\nu(Q)\le\nnti .
\]
Every point lies on $5$ blocks, so $6\nnti\le\sum_B|B\cap\X|=5x$ (Lemma~\ref{lem:blocksid}).

If $a_0\le2$, then $s-a_0\ge61-x$ and hence $2(61-x)\le\frac56x$, i.e.\ $17x\ge732$, and $x\ge44$.

Let $a_0\ge3$. By Lemma~\ref{lem:A0}, $\A'_0$ lies in a plane $W$. All blocks through points of $\A'_0$ are TI, and we count them. For a block $B$ let $k_B=|B\cap\A'_0|$. A block with $k_B\ge1$ is counted $k_B$ times in $5a_0$, so the number of blocks meeting $\A'_0$ is $5a_0-c$ with the excess $c=\sum_B(k_B-1)^+$.

\emph{Case $W$ is a block.} Every other block meets $W$ in at most a point, so $c=a_0-1$ and there are $4a_0+1\ge2a_0+6$ blocks meeting $\A'_0$.

\emph{Case $W$ is not a block.} A block $B$ with $k_B\ge2$ meets $W$ in the line $L_B$ spanned by two points of $\A'_0$, and $k_B=|L_B\cap\A'_0|$. Distinct blocks give distinct lines, since a line lies in at most one block. For a line $L$ of $W$ put $k_L=|L\cap\A'_0|$ and let $z$ be the number of lines of $W$ disjoint from $\A'_0$. Each point of $W$ lies on $3$ of its $7$ lines, so $\sum_Lk_L=3a_0$ and
\[
c\le\sum_{L\le W}(k_L-1)^+=\sum_{L\le W}(k_L-1)+z=3a_0-7+z .
\]
Two lines of $W$ cover $5$ points and $|W\setminus\A'_0|\le4$, so $z\le1$. So at least $5a_0-(3a_0-7+z)\ge2a_0+6$ blocks meet $\A'_0$, and all of them are TI.

In both cases $n_{\mathrm{TI}}\ge2a_0+6$. With $\nnti\ge2(s-a_0)$ this gives
\[
45=n_{\mathrm{TI}}+\nnti\ge2s+6,
\]
so $s\le19$ and $x=63-s\ge44$.
\end{proof}

\section{\texorpdfstring{The bound $|\X_H|\ge41$}{The bound |XH| >= 41}}\label{sec:41}

The polarity of Section~\ref{sec:polar} is available as long as $\A'$ spans $H$ and $|\X_H|\le47$ (Lemma~\ref{lem:polspan}). In this section we treat the case that $\A'$ does not span $H$ and then combine both cases.

Suppose that $\A'$ lies in a $5$-subspace $F$. For $|\X_H|\le40$, $F$ is a $\beta$-flat (Lemma~\ref{lem:nonspan}), and we show that this is impossible. The $40$ blocks outside $F$ meet $F$ in lines, their traces, and these form a rigid family of lines (Lemma~\ref{lem:cosets}). A count of the points outside $F$ then gives a contradiction (Lemmas~\ref{lem:Opairs} and \ref{lem:ea}, Proposition~\ref{prop:caseN}).

\begin{lemma}\label{lem:formF}
Suppose that $\A'$ spans a $5$-subspace $F$ of $H$ and $|\A'|\ge8$. Then there is an alternating form on $F$ with $Q^\perp=G(Q)\cap F$ for every $Q\in\A'$ \textup{(}orthogonality inside $F$\textup{)}, and its radical is a point $r$. If $F$ is a $\beta$-flat, then $r$ is its focal point $Q_0\in\A'$. Otherwise $r\notin\A'$. For $Q\in\A'$ with $G(Q)\ne F$, $Q^\perp$ is a solid.
\end{lemma}
\begin{proof}
As in the proof of Lemma~\ref{lem:polspan}, now with the restrictions $a_Q|_F$ (with kernel $G(Q)\cap F$), $\langle u,v\rangle=M(u)(v)$ defines an alternating form on $F$ with $Q^\perp=G(Q)\cap F$ for $Q\in\A'$. If $G(Q)\ne F$, this is a solid, since two distinct $5$-subspaces of $H$ meet in a solid. A point $Q\in\A'$ lies in the radical $R$ iff $Q^\perp=F$ iff $G(Q)=F$, and this holds for at most one $Q$ by Lemma~\ref{lem:dh}(b), namely the focal point $Q_0$ if $F$ is a $\beta$-flat. As $\dim F=5$ and the rank is even, $\dim R\in\{1,3,5\}$. If $\dim R=5$, then $G(Q)=F$ for all $Q\in\A'$, so $|\A'|\le1$. If $\dim R=3$, the $24$ points of $F\setminus R$ fall into $3$ cosets of $R$. Points of the same coset have the same solid $Q^\perp$, and only two $5$-subspaces of $H$ other than $F$ contain a given solid of $F$, so each coset contains at most $2$ points of $\A'$ by Lemma~\ref{lem:dh}(b). Then $|\A'|\le1+3\cdot2=7$. Both cases contradict $|\A'|\ge8$, so $R$ is a point $r$.
\end{proof}

\begin{lemma}\label{lem:nonspan}
Suppose that $\A'$ is contained in a $5$-subspace $F$ of $H$. Then every $Q\in\A'$ with $G(Q)\ne F$ lies on exactly one block of $F$, and this block lies in $G(Q)\cap F$. This gives:
\begin{enumerate}[label=(\alph*)]
\item if $F$ is not a $\beta$-flat, then $|\A'|\le21$, i.e.\ $|\X_H|\ge42$;
\item if $\A'$ spans a subspace of dimension at most $4$, then $|\A'|\le15$, i.e.\ $|\X_H|\ge48$;
\item if $F=\langle\A'\rangle$ and $|\X_H|\le40$, then $F=G(Q_0)$ is a $\beta$-flat whose focal point $Q_0$ lies in $\A'$, and $|\X_H|=32+t$ with $t=|\X\cap F|\le8$.
\end{enumerate}
\end{lemma}
\begin{proof}
Let $Q\in\A'$ with $G(Q)\ne F$. Then $G(Q)\cap F$ is a solid of the $\beta$-flat $G(Q)$ through its focal point $Q$, hence rich. Its block lies in $G(Q)$ and therefore passes through $Q$. By Lemma~\ref{lem:focal}, $j_F(Q)=1$.

(a) If $F$ is not a $\beta$-flat, every $Q\in\A'$ has $j_F(Q)=1$, so $|\A'|\le m_1(F)\le23$ by Proposition~\ref{prop:patterns}, and $m_1(F)\le21$ unless $F$ has pattern $t_2$. Suppose that $|\A'|\ge22$. Then $\A'$ spans $F$, since a subspace of dimension at most $4$ has at most $15$ points. By Lemma~\ref{lem:formF} there is an alternating form on $F$ with $Q^\perp=G(Q)\cap F$ for $Q\in\A'$ whose radical is a point $r\notin\A'$, since $F$ is not a $\beta$-flat. If a subspace $U$ of $F$ is totally isotropic, so is $U+\langle r\rangle$, which has dimension at most $3$. So every totally isotropic plane of $F$ contains $r$. For a block $B\le F$ and $Q\in B\cap\A'$ we have $B\le G(Q)\cap F=Q^\perp$, so $B\cap\A'$ lies in the radical $B\cap B^\perp$ of the restriction of the form to $B$.

Since $|\A'|\ge22>21$, $F$ has pattern $t_2$. Four blocks $B_1,\dots,B_4$ of $F$ pass through the point $P_4$ with $j_F(P_4)=4$ and the fifth block $B_5$ avoids $P_4$. Each of $B_1,\dots,B_4$ contains $5$ points with $j_F=1$ and $B_5$ contains $3$. This gives the $23$ points with $j_F=1$, and at most one of them is not in $\A'$. So $B_1$ and $B_2$ contain at least $4$ points of $\A'$, which span them, and therefore $B_1$ and $B_2$ are totally isotropic. Moreover $|B_5\cap\A'|\ge2$, so $B_5\cap B_5^\perp$ has dimension at least $2$, and $B_5$ is totally isotropic since a plane has rank $0$ or $2$. Hence $r\in B_1\cap B_2\cap B_5=\{P_4\}\cap B_5=\emptyset$, a contradiction. So $|\A'|\le21$.

(b) A subspace of dimension at most $4$ has at most $15$ points.

(c) As $|\X_H|\le40<42$, $F$ is a $\beta$-flat by (a). Its focal point is the focal point of a $\beta$-flat in $H$, so it lies in $\A'$. Finally $\A'=F\setminus\X$ has $31-t$ points, so $|\X_H|=63-(31-t)$.
\end{proof}

\begin{definition}
The \emph{flat case} is: $|\X_H|\le40$ and $\A'$ does not span $H$.
\end{definition}

In the flat case $\A'$ lies in a proper subspace and $|\A'|=63-|\X_H|\ge23>15$, so $F=\langle\A'\rangle$ is a $5$-subspace. By Lemma~\ref{lem:nonspan}(c), $F=G(Q_0)$ is a $\beta$-flat with $Q_0\in\A'$, and $|\X_H|=32+t$ with $t=|\X\cap F|\le8$. Until Proposition~\ref{prop:caseN} we assume the flat case, with the aim of a contradiction. We call the $5$ blocks of $F$ the \emph{inner} blocks and the other $40$ blocks the \emph{outer} blocks. An outer block meets $F$ in a line, its \emph{trace}.

\begin{lemma}\label{lem:caseNform}
The following hold.
\begin{enumerate}[label=(\alph*)]
\item There is an alternating form on $F$ such that $Q^\perp=G(Q)\cap F$ \textup{(}a solid\textup{)} for every $Q\in\A'\setminus\{Q_0\}$, and its radical is $R=\langle Q_0\rangle$, so the form induces a nondegenerate form on $F/Q_0\cong\PG(3,2)$.
\item The images of the inner blocks form a line spread $S_0$ of $F/Q_0$, and $S_0$ is a spread of $W(2)$.
\item Distinct outer blocks have distinct traces, a trace does not pass through $Q_0$ and is not contained in an inner block, and every point of $F\setminus\{Q_0\}$ lies on exactly $4$ traces.
\item Let $\Lambda$ be the set of lines of $F$ not through $Q_0$ whose image in $F/Q_0$ is a TI line not in $S_0$. Then $|\Lambda|=40$, and every point of $F\setminus\{Q_0\}$ lies on exactly $4$ lines of $\Lambda$. For $Q\in\A'\setminus\{Q_0\}$ the traces through $Q$ are exactly the $4$ lines of $\Lambda$ through $Q$. Every trace which contains no point of $\A'$ lies in $\X\cap F$.
\end{enumerate}
\end{lemma}
\begin{proof}
(a) $\A'$ spans $F$ and $|\A'|\ge23$. By Lemma~\ref{lem:formF} there is an alternating form on $F$ with $Q^\perp=G(Q)\cap F$ for $Q\in\A'$, and its radical is the focal point $Q_0$ of the $\beta$-flat $F=G(Q_0)$. For $Q\in\A'\setminus\{Q_0\}$ we have $G(Q)\ne F$, so $Q^\perp$ is a solid.

(b) By Lemma~\ref{lem:dh}(c), the inner blocks pass through $Q_0$ and partition $F\setminus\{Q_0\}$, so $S_0$ is a line spread. An inner block $B$ containing a point $Q\in\A'\setminus\{Q_0\}$ lies in $G(Q)$ (Lemma~\ref{lem:GQ}), hence in $Q^\perp$, so its image is TI. Since $t=|\X_H|-32\le8<12$, at most one inner block has all its $6$ points other than $Q_0$ in $\X$, so at least $4$ lines of $S_0$ are TI. Suppose that a line $m\in S_0$ is not TI, and let $p\in m$. The plane $p^\perp$ of $F/Q_0$ contains exactly one line $\ell$ of $S_0$ (Lemma~\ref{lem:spreadplane}), and $\ell\ne m$, since $m$ contains a point $q$ with $\langle p,q\rangle=1$. If $\ell$ were TI, then $\ell^\perp=\ell$, and $\ell\le p^\perp$ would give $p\in\ell$, which is false. So a non-TI line of $S_0$ forces a second one, contradicting the fact that at most one line of $S_0$ is non-TI. Therefore all $5$ lines of $S_0$ are TI.

(c) An outer block $B\not\le F$ meets $F$ in a line. The blocks through $Q_0$ are the inner blocks (Lemma~\ref{lem:dh}(c)), and a line lies in at most one block. A point $P\in F\setminus\{Q_0\}$ lies on $5$ blocks, exactly one of them inner. The traces of the other $4$ are distinct lines through $P$, and every trace through $P$ arises in this way.

(d) There are $10$ TI lines not in $S_0$. The preimage of such a line $L$ is the plane $\langle L,Q_0\rangle$, whose $4$ lines not through $Q_0$ map onto $L$. A point $P\ne Q_0$ has an image on $3$ TI lines, one of them in $S_0$. For each of the other two, $P$ lies on $2$ lines of $\langle L,Q_0\rangle$ not through $Q_0$. Let $Q\in\A'\setminus\{Q_0\}$. Its $4$ outer blocks lie in $G(Q)$, so their traces are lines through $Q$ in $G(Q)\cap F=Q^\perp$, not through $Q_0$. Such a line is TI, and its image is not in $S_0$ (otherwise the trace would lie in an inner block). Therefore the $4$ traces through $Q$ lie in $\Lambda$, which has exactly $4$ lines through $Q$. The last statement holds because $F\setminus\{Q_0\}=(\A'\setminus\{Q_0\})\cup(\X\cap F)$.
\end{proof}

The incidence structure behind $\Lambda$ is summarized in the following table. Here $F/Q_0\cong\PG(3,2)$ carries the nondegenerate form of Lemma~\ref{lem:caseNform}(a), and $\pi_L=\langle L,Q_0\rangle$ denotes the plane of $F$ over a line $L$ of $F/Q_0$.
\begin{center}\small
\begin{tabular}{@{}lll@{}}
 & in $F/Q_0$ & in $F$\\\hline
lines of $S_0$ & $5$ TI lines & the inner blocks, through $Q_0$\\
other TI lines $L$ & $10$ lines & $10$ planes $\pi_L$, each with $4$ lines of $\Lambda$ \textup{(}$|\Lambda|=40$\textup{)}\\
through a point $\bar P$ & $3$ TI lines, $1$ in $S_0$ & $P$ lies on $2+2=4$ lines of $\Lambda$, two in each of two planes $\pi_L$
\end{tabular}
\end{center}

\begin{lemma}\label{lem:cosets}
The $40$ traces are exactly the $40$ lines of $\Lambda$. Consequently, for every $P\in F\setminus\{Q_0\}$ the following holds. Let $L_1,L_2$ be the two TI lines of $F/Q_0$ through the image of $P$ which are not in $S_0$, and let $\pi_i=\langle L_i,Q_0\rangle$. Then $\pi_1\cap\pi_2$ is the line through $P$ and $Q_0$, and the four traces through $P$ are the two lines through $P$ not through $Q_0$ in $\pi_1$ and the two such lines in $\pi_2$.
\end{lemma}
\begin{proof}
Let $U$ be the set of lines of $\Lambda$ which are not traces and $W$ the set of traces not in $\Lambda$. Then $U\cap W=\emptyset$. No line of $U\cup W$ passes through $Q_0$ (by the definition of $\Lambda$ and Lemma~\ref{lem:caseNform}(c)), and by Lemma~\ref{lem:caseNform}(d) none contains a point of $\A'\setminus\{Q_0\}$. So all lines of $U\cup W$ lie in the $t$-set $\X\cap F$. For $P\in F\setminus\{Q_0\}$ let $u(P)$ and $w(P)$ be the numbers of lines of $U$ and of $W$ through $P$. Since $P$ lies on exactly $4$ traces and on exactly $4$ lines of $\Lambda$, both numbers equal $4$ minus the number of traces in $\Lambda$ through $P$, so $u(P)=w(P)$. If $u(P)\ge2$, then $P$ lies on at least $4$ distinct lines of $U\cup W$, which pairwise meet only in $P$ and therefore cover at least $1+4\cdot2=9>t$ points of $\X\cap F$. So $u(P)=w(P)\le1$ for every point $P$. Suppose that $U$ contains a line $u$. Each of the three points of $u$ lies on exactly one line of $W$, and these three lines are distinct, since a line of $W$ through two points of $u$ would equal $u$. As $w(P)\le1$ for all $P$, they are pairwise disjoint and cover $9>t$ points of $\X\cap F$, a contradiction. Thus $U=\emptyset$, so $w(P)=u(P)=0$ for all $P$, and $W=\emptyset$.

The second statement is the description of the lines of $\Lambda$ through $P$ in the proof of Lemma~\ref{lem:caseNform}(d). The planes $\pi_1\ne\pi_2$ both contain the line through $P$ and $Q_0$.
\end{proof}

\begin{definition}
For $P\in F\setminus\{Q_0\}$ the five blocks through $P$ pairwise meet only in $P$, and one of them is inner. Let $O(P)$ be the set of the $16$ points of $H\setminus F$ which lie on the four outer blocks through $P$. For $a\in H\setminus F$ let $T(a)$ be the union of the traces of the $5$ blocks through $a$, which are all outer.
\end{definition}

Thus $T(a)$ is the shadow of $a$ on $F$, and $O(P)$ is the set of points $a$ whose shadow contains $P$. The $5$ blocks through $a$ pairwise meet only in $a\notin F$, and no trace passes through $Q_0$ (Lemma~\ref{lem:caseNform}(c)), so $T(a)\subseteq F\setminus\{Q_0\}$ is a disjoint union of $5$ lines and $|T(a)|=15$. The $30$ points of $F\setminus\{Q_0\}$ fall into the $15$ pairs $\{P,P+Q_0\}$.

\begin{lemma}\label{lem:Opairs}
Let $P\in F\setminus\{Q_0\}$ and $P'=P+Q_0$.
\begin{enumerate}[label=(\alph*)]
\item For $a\in H\setminus F$ we have $P\in T(a)$ iff $a\in O(P)$. If $P\in\A'$, then $O(P)=G(P)\setminus F$.
\item If $P\in\X$ and $P'\in\A'$, then $O(P)=(H\setminus F)\setminus G(P')$.
\item If $\{P,P'\}$ is not contained in $\X$, then every set $T(a)$, $a\in H\setminus F$, contains exactly one of $P,P'$.
\item $|O(P)\cap O(P')|\le8$.
\end{enumerate}
\end{lemma}
\begin{proof}
(a) A block through $a\notin F$ is outer, so $P\in T(a)$ iff $P$ and $a$ lie on a common outer block, i.e.\ iff $a\in O(P)$. If $P\in\A'$, then $G(P)\ne F$ since $P\ne Q_0$, so $G(P)\cap F$ is a solid and $|G(P)\setminus F|=16$. The blocks through $P$ lie in $G(P)$ (Lemma~\ref{lem:GQ}), hence $O(P)\subseteq G(P)\setminus F$, and both sets have $16$ points.

(b) The form of Lemma~\ref{lem:caseNform}(a) has radical $\langle Q_0\rangle$, so $P^\perp=P'^\perp=G(P')\cap F$, a solid of $F$ through $Q_0$. Let $B$ be an outer block through $P$ with trace $\ell$. By Lemma~\ref{lem:cosets}, $\ell\in\Lambda$, so the image of $\ell$ is a TI line through the image of $P$ and $\ell\le P^\perp=G(P')\cap F$. If $B\le G(P')$, then $B$ is a block of the $\beta$-flat $G(P')$ and passes through $P'$ (Lemma~\ref{lem:dh}(c)), so $B$ contains $P+P'=Q_0$ and is an inner block, which is false. Hence $B\cap G(P')=\ell$, and $B$ has no point in $G(P')\setminus F$. So $O(P)\subseteq(H\setminus F)\setminus G(P')$, and both sets have $16$ points.

(c) If $P,P'\in\A'$, then $G(P)\ne G(P')$ (Lemma~\ref{lem:dh}(b)) are $5$-subspaces different from $F$ through the solid $P^\perp=P'^\perp$. With $F$ they are the three $5$-subspaces of $H$ through this solid, so $G(P)\setminus F$ and $G(P')\setminus F$ partition $H\setminus F$, and (a) gives the claim. If $P\in\X$ and $P'\in\A'$, then by (a) and (b), $P\in T(a)$ iff $a\in O(P)$ iff $a\notin G(P')$ iff $P'\notin T(a)$. The case $P\in\A'$, $P'\in\X$ is symmetric.

(d) Let $\pi_1,\pi_2$ be as in Lemma~\ref{lem:cosets}. It applies to $P$ and to $P'$ with the same planes, since $P$ and $P'$ have the same image in $F/Q_0$. Let $B\ni P$ and $B'\ni P'$ be outer blocks with traces $\ell$ and $\ell'$. Then $\ell\ne\ell'$, as no trace contains $P$, $P'$ and hence $Q_0$, and $B\ne B'$, as otherwise $B$ would contain $\ell+\ell'$ and hence $Q_0$. So $B\cap B'$ is at most a point. If $\ell$ and $\ell'$ lie in the same plane $\pi_i$, they meet in a point of $F$, and $B$ and $B'$ have no common point outside $F$. Each of $P,P'$ has two traces in each $\pi_i$, so there are $2\cdot2+2\cdot2=8$ pairs $(B,B')$ with $\ell$ and $\ell'$ in different planes, and each of them contributes at most one point to $O(P)\cap O(P')$.
\end{proof}

\begin{lemma}\label{lem:ea}
For $a\in H\setminus F$ let $e_a$ be the number of pairs $\{P,P+Q_0\}$ contained in $T(a)$. Then exactly $e_a$ pairs are disjoint from $T(a)$, all $2e_a$ of these pairs lie in $\X\cap F$, and $e_a\ge2$.
\end{lemma}
\begin{proof}
Let $c$ be the number of pairs with exactly one point in $T(a)$. Then $2e_a+c=|T(a)|=15$ and the number of pairs disjoint from $T(a)$ is $15-e_a-c=e_a$. A pair which is not contained in $\X$ meets $T(a)$ in exactly one point (Lemma~\ref{lem:Opairs}(c)), so all pairs of the other two kinds lie in $\X\cap F$.

The $5$ traces through $a$ are pairwise disjoint lines of $\Lambda$ (Lemma~\ref{lem:cosets}). Two of them cannot have the same image $L$ in $F/Q_0$, since both would lie in the plane $\langle L,Q_0\rangle$ and meet. So their images are $5$ distinct totally isotropic lines not in $S_0$, and no point of $F/Q_0$ lies on three of them, since only two TI lines outside $S_0$ pass through a point. A point of $F/Q_0$ has both preimages in $T(a)$ iff it lies on two image lines, so $e_a$ is the number of pairs of meeting image lines, and the images cover $15-e_a$ points. If $e_a=0$, the images form a spread of $W(2)$ disjoint from the spread $S_0$ (Lemma~\ref{lem:caseNform}(b)), contradicting Lemma~\ref{lem:W2}. If $e_a=1$, the images cover all points of $F/Q_0$ but one point $r$. A TI line $L\not\ni r$ meets the plane $r^\perp$ in exactly one point, since $L\le r^\perp$ would give $r\in L^\perp=L$. So the $5$ images cover at most $5$ of the $6$ points of $r^\perp\setminus\{r\}$, a contradiction.
\end{proof}

\begin{proposition}\label{prop:caseN}
The flat case does not occur, i.e.\ if $32\le|\X_H|\le40$, then $\A'$ spans $H$.
\end{proposition}
\begin{proof}
Let $k$ be the number of pairs $\{P,P+Q_0\}$ contained in $\X\cap F$. Then $k\le\lfloor t/2\rfloor\le4$, since $t\le8$ (Lemma~\ref{lem:nonspan}). For a pair $\{P,P'\}\subseteq\X\cap F$ and $a\in H\setminus F$ we have $\{P,P'\}\cap T(a)=\emptyset$ iff $a\notin O(P)\cup O(P')$ (Lemma~\ref{lem:Opairs}(a)). By Lemma~\ref{lem:Opairs}(d), $|O(P)\cup O(P')|\ge16+16-8=24$, so this holds for at most $8$ of the $32$ points $a$. Count the pairs $(a,\{P,P'\})$ with $a\in H\setminus F$, $\{P,P'\}\subseteq\X\cap F$ and $\{P,P'\}\cap T(a)=\emptyset$. By Lemma~\ref{lem:ea} there are at least $\sum_ae_a\ge64$ of them, and by the above at most $8k\le32$. This is a contradiction.
\end{proof}

\begin{theorem}\label{thm:41}
For every hyperplane section, $|\X_H|=|\Y_H|\ge41$, so a hyperplane section has at most $22$ $\beta$-flats. If $\A'$ spans $H$, then even $|\X_H|\ge44$.
\end{theorem}
\begin{proof}
If $\A'$ spans $H$, then either $|\X_H|\ge48$ or Lemma~\ref{lem:polspan} and Theorem~\ref{thm:pol44} give $|\X_H|\ge44$. Otherwise $\A'$ lies in a proper subspace, which has at most $31$ points, so $|\X_H|\ge32$, and Proposition~\ref{prop:caseN} excludes $32\le|\X_H|\le40$.
\end{proof}

\begin{lemma}\label{lem:griesmer}
Let $\N$ be a set of points of $\PG(6,2)$ with $|\N\cap H|\ge t\ge1$ for every hyperplane $H$. Then $|\N|\ge n(t)$, the least $n$ with $127-n\ge g(64-n+t)$, where $g(d)=\sum_{i=0}^{6}\lceil d/2^i\rceil$. In particular $n(41)=85$.
\end{lemma}
\begin{proof}
Let $n=|\N|$ and let $M$ be the $7\times(127-n)$ matrix whose columns are the points outside $\N$. For $0\ne u\in\F_2^7$ the codeword $uM$ has weight $\#\{P\notin\N: u\cdot P\ne0\}=64-|\N\setminus u^\perp|=64-n+|\N\cap u^\perp|\ge 64-n+t$. If $64-n+t\ge1$, $M$ generates a binary linear $[127-n,7,\ge64-n+t]$ code, and the Griesmer bound \cite{Gri60} gives $127-n\ge g(64-n+t)$. As $g(d)-g(d-1)\ge1$, the function $n\mapsto127-n-g(64-n+t)$ is non-decreasing, so the inequality holds exactly for $n\ge n(t)$. For $t=41$ and $n=84$ the condition reads $43\ge g(21)=45$, which fails, and for $n=85$ it reads $42\ge g(20)=42$, which holds.
\end{proof}

\begin{corollary}\label{cor:41}
Let $\D$ be a Fano design. Every hyperplane contains at least $41$ non-$\alpha$-points, hence at most $22$ $\alpha$-points. Moreover, $|\N|\ge85$ and $|\A|\le42$, and at least $\lceil127\cdot41/3\rceil=1736$ of the $2667$ five-subspaces are not $\beta$-flats.
\end{corollary}
\begin{proof}
Theorem~\ref{thm:41}, Lemma~\ref{lem:K}(c) and Lemma~\ref{lem:griesmer} with $t=41$. Every $5$-subspace lies in $3$ hyperplanes, and $\sum_H|\Y_H|\ge127\cdot41$.
\end{proof}

\section{Further results on hyperplane sections}\label{sec:dh}

Hyperplane sections exist (Observation~\ref{obs:dh}), and the local conditions at a single point or $5$-subspace impose no restriction.

\subsection{Existence and local conditions}

\begin{observation}\label{obs:dh}
Hyperplane sections exist, and an explicit example is given in Appendix~\ref{app:section}. We know $60$ pairwise distinct hyperplane sections, found by integer programming with prescribed nontrivial automorphism groups and by random sampling of such programs (we did not test them for isomorphism). A second program, written independently of the search code, confirms that all $60$ are hyperplane sections and that their numbers of $\beta$-flats are $0$ ($42$ sections), $1$ ($5$), $2$ ($2$), $3$ ($2$), $6$ ($1$) and $7$ ($8$), so $|\X_H|\ge56$ for all of them, far from the bound $|\X_H|\ge41$. A further section with $15$ $\beta$-flats and $|\X_H|=48$ is given in Appendix~\ref{app:section}. Its focal points form a solid, so Lemma~\ref{lem:nonspan}(b) is sharp.
\end{observation}

An automorphism of $\D|_H$ need not extend to $\D$, so these symmetric configurations do not contradict \cite{KKW18}. However, by \cite{KKW18} a lift of such a configuration to a $q$-analog of the Fano plane cannot be invariant under a group of order larger than $2$.

In a design, the $5$ blocks of $\D|_H$ through a point $P\in H$ give $5$ pairwise disjoint lines of $H/P\cong\PG(4,2)$, which belong to the line spread $\Der_P(\D)$ of $V/P\cong\PG(5,2)$.

\begin{fact}\label{fact:local}
Every set of $5$ pairwise disjoint lines of a hyperplane $W\cong\PG(4,2)$ of $\PG(5,2)$ extends to a line spread of $\PG(5,2)$, and each of the $10$ orbits occurs as the set of spread lines inside a hyperplane for some of the $131\,044$ isomorphism types of line spreads of $\PG(5,2)$ \cite{MT09}.
\end{fact}

Together with Lemma~\ref{lem:dh}(a) (every $5$-subspace automatically contains exactly $5$ blocks), the local conditions at points and $5$-subspaces of $H$ are therefore automatic, and any obstruction cannot come from a single point or a single $5$-subspace. The structure around a single point gives no obstruction either: Etzion and Hooker \cite{EH18} construct, for every $q$, a derived and a residual design with the parameters coming from a $q$-analog of the Fano plane at a point, and ask whether such a structure can be extended to a design. For a hyperplane section the analogous extension question is our lifting problem (Proposition~\ref{prop:lift}).

\subsection{Lifting a hyperplane section}
We now ask when a hyperplane section extends to a Fano design. Proposition~\ref{prop:lift} reduces this to the existence of a Steiner system $S(2,4,64)$ \textup{(}compare \cite[Thm.~3.1]{EV11}\textup{)} on the complement of $H$ whose directions are the free lines.

Let $H$ be a hyperplane of $V=\F_2^7$ and $\K$ a set of planes of $H$. Identify $V\setminus H$ with the affine space $\AG(6,2)$. For a line $L$ of $H$ write $\bar L=L\cup\{0\}$ for the corresponding $2$-dimensional subspace.

\begin{proposition}\label{prop:lift}
A set $\D$ of planes of $V$ with $\D|_H=\K$ is a $2$-$(7,3,1)_2$ design iff $\K$ is a hyperplane section and the blocks of $\D$ not in $H$ are the planes $\langle a,L\rangle=\bar L\cup(a+\bar L)$, where the affine planes $a+\bar L$ form a Steiner system $S(2,4,64)$ on $\AG(6,2)$ whose directions $L$ are exactly the free lines of $\K$, each occurring once.
\end{proposition}

\begin{proof}
A block $B\not\le H$ meets $H$ in a line $L$, and $B\setminus H=a+\bar L$ is an affine plane with direction $L$. It contains exactly one line of $H$ (namely $L$) and the $6$ lines $\{b,b',b+b'\}$, $b\ne b'\in a+\bar L$, which correspond to the pairs of affine points. Every line of $V$ either lies in $H$ or is of the form $\{b,b',b+b'\}$ with $b,b'\notin H$.

Hence $\D$ is a design iff (i) every line of $H$ lies in exactly one block of $\K$ or is the direction of exactly one lifted block, and (ii) every pair of affine points lies in exactly one lifted block. Given (ii), the pairs $\{b,b+d\}$ with difference $d\in H$ ($32$ of them) are covered exactly once, and a lifted block covers $2$ of them if $d$ lies on its direction and none otherwise, so exactly $16$ directions pass through $d$. With (i), $d$ lies on $(31-16)/3=5$ blocks of $\K$, and no line of $H$ lies in two blocks of $\K$: $\K$ is a hyperplane section.

Conversely, if $\K$ is a hyperplane section, (i) says that the directions are exactly the free lines, each once, and (ii) is the Steiner system condition.
\end{proof}

\begin{corollary}
A lift of a hyperplane section $\K$ is the same as a coloring $c$ of the $336$ free lines of $\K$, where $c(L)$ is one of the $16$ cosets of $\bar L$ in $V\setminus H$, such that for every point $d$ of $H$ the $16$ free lines $L$ through $d$ give $32$ pairwise different cosets of $\langle d\rangle$ \textup{(}each $c(L)$ is the union of two cosets of $\langle d\rangle$\textup{)}.
\end{corollary}
\begin{proof}
By Proposition~\ref{prop:lift}, a lift chooses, for each free line $L$, one of the $16$ affine planes $a+\bar L$ with direction $L$. The coloring condition is (ii) sorted by differences: the pairs with difference $d$ are the $32$ cosets of $\langle d\rangle$ in $V\setminus H$, and the lifted block with direction $L\ni d$ contains exactly two of them.
\end{proof}

Counting gives $|\K|=45$ automatically ($7|\K|+336=651$). At an affine point $a$, the directions of the $21$ lifted blocks through $a$ form a line spread of $H$ consisting of free lines.

\subsection{\texorpdfstring{Three hyperplane sections through a $5$-subspace}{Three hyperplane sections through a 5-subspace}}
Let $\D$ be a Fano design, $F$ a $5$-subspace, $H_1,H_2,H_3$ the three hyperplanes containing $F$, $j(P)=j_F(P)$, and $T_i=\{B\cap F: B\in\D,\ B\le H_i,\ B\not\le F\}$, the \emph{trace} of $\D|_{H_i}$ on $F$. In coordinates, the lines disjoint from $F$ become the cells of a $32\times32$ Latin square, and each transversal block occupies one $2\times2$ subsquare (Lemma~\ref{lem:pencil}(d)).

\begin{lemma}\label{lem:pencil}
\begin{enumerate}[label=(\alph*)]
\item Every block meets $F$ in a point, a line or a plane. Exactly $5$ blocks lie in $F$, exactly $120$ meet $F$ in a line and exactly $256$ meet $F$ in a point. A block meeting $F$ in a line lies in exactly one $H_i$, and each $H_i$ contains $40$ of them. A block meeting $F$ in a point (\emph{transversal block}) lies in no $H_i$.
\item The $155$ lines of $F$ are partitioned into the $35$ lines of the blocks in $F$ and $T_1,T_2,T_3$ ($40$ lines each). Every point $P\in F$ lies on exactly $5-j(P)$ lines of each $T_i$.
\item Through $P\in F$ there are exactly $6+2j(P)$ transversal blocks, and $H_i$ contains exactly $6+2j(P)$ free lines through $P$ not in $F$. A transversal block $B$ through $P$ meets each $H_i$ in a free line $L_i\ni P$, and $L_1,L_2,L_3$ are the three lines of $B$ through $P$.
\item The $1024$ lines disjoint from $F$ are covered by the transversal blocks, $4$ by each. Write $V=F\oplus\langle e_1,e_2\rangle$ with $H_1=F+e_1$, $H_2=F+e_2$. A line disjoint from $F$ is $\{e_1+x,e_2+y,e_1+e_2+x+y\}$ for a unique cell $(x,y)\in F\times F$ of the Cayley table of $(F,+)\cong\F_2^5$, and the $4$ lines of a transversal block through $P$ which do not pass through $P$ form an intercalate $\{x,x+P\}\times\{y,y+P\}$.
\end{enumerate}
\end{lemma}
\begin{proof}
(a) Since $\dim(B\cap F)\ge3+5-7$, every block meets $F$. The $5$ blocks in $F$ contain $35$ of the $155$ lines of $F$, and each of the other $120$ lines lies in a unique block, which contains no other line of $F$. Counting incidences of the blocks with the $31$ points of $F$, each of which lies on $21$ blocks, gives $31\cdot21=5\cdot7+120\cdot3+256\cdot1$. A block meeting $F$ in a line spans a hyperplane with $F$, and a block meeting $F$ in a point spans $V$. Each $H_i$ contains $45$ blocks, $5$ of them in $F$.

(b) The traces of different blocks are different lines, and a block lies in only one $H_i$. $P$ lies on $5$ blocks of $H_i$, $j(P)$ of them in $F$, and each of the other $5-j(P)$ meets $F$ in a line through $P$.

(c) Of the $21$ blocks through $P$, $j(P)$ lie in $F$ and $3(5-j(P))$ meet $F$ in a line, so $6+2j(P)$ are transversal. In $H_i$, $P$ lies on $16$ lines not in $F$. The $5-j(P)$ blocks of $H_i$ through $P$ not in $F$ contain $2$ of them each, so $6+2j(P)$ are free. $B\cap H_i$ is a line through $P=B\cap F$. It is free in $H_i$ since $B\not\le H_i$, and $L_i\cap L_k\le B\cap F=P$.

(d) Let $\ell$ be a line disjoint from $F$. Its image in $V/F$ consists of the three nonzero vectors, so $\ell$ has one point in each $H_i\setminus F$. Let $B$ be the block through $\ell$. Any two lines of $B$ meet and $\ell$ is disjoint from $F$, so $B\cap F$ is not a line. Since $B\not\le F$, it follows that $B\cap F$ is a single point $P$, so $B$ is transversal. Conversely, a transversal block $B$ with $B\cap F=\{P\}$ has $7$ lines, $3$ through $P$ and $4$ not through $P$, and the latter are disjoint from $F$. The number of lines disjoint from $F$ is $2667-155-31\cdot48=1024$, since $155$ lines lie in $F$ and each of the $31$ points of $F$ lies on $48$ lines not in $F$. There are $256$ transversal blocks, and $1024=4\cdot256$.

For the coordinates, write $B=\langle P,e_1+x,e_2+y\rangle$ with $x,y\in F$. Its points outside $F$ are $e_1+x+\{0,P\}$, $e_2+y+\{0,P\}$ and $e_1+e_2+x+y+\{0,P\}$. Each of its $4$ lines not through $P$ has the form $\{e_1+x',e_2+y',e_1+e_2+x'+y'\}$ with $x'\in\{x,x+P\}$ and $y'\in\{y,y+P\}$, that is, it joins a point of $\{e_1+x,e_1+x+P\}$ to a point of $\{e_2+y,e_2+y+P\}$. This is the cell $(x',y')$, and the $4$ cells form the intercalate $\{x,x+P\}\times\{y,y+P\}$. For example, for $x=y=0$ the cells $(0,0)$, $(0,P)$, $(P,0)$, $(P,P)$ are the lines $\{e_1,e_2,e_1+e_2\}$, $\{e_1,e_2+P,e_1+e_2+P\}$, $\{e_1+P,e_2,e_1+e_2+P\}$ and $\{e_1+P,e_2+P,e_1+e_2\}$.
\end{proof}

\begin{proposition}\label{prop:pencil}
$\D=\D|_{H_1}\cup\D|_{H_2}\cup\D|_{H_3}\cup\{\text{transversal blocks}\}$. Conversely, let $\K_i$ be hyperplane sections on $H_i$ \textup{(}$i=1,2,3$\textup{)} with the same $5$ planes in $F$ and with $T_1,T_2,T_3$ pairwise disjoint (then they partition the $120$ lines of $F$ not in these $5$ planes). Then $\K_1\cup\K_2\cup\K_3$ extends to a Fano design iff the $1024$ cells of the Cayley table can be partitioned into $256$ intercalates $\{x,x+P\}\times\{y,y+P\}$, $0\ne P\in F$, such that every free line of $H_1$, $H_2$, $H_3$ outside $F$ is used exactly once. Here the intercalate $\{x,x+P\}\times\{y,y+P\}$ uses the three lines
\[
\{P,e_1+x,e_1+x+P\},\quad \{P,e_2+y,e_2+y+P\},\quad \{P,e_1+e_2+x+y,e_1+e_2+x+y+P\}
\]
of $H_1$, $H_2$, $H_3$, respectively.
\end{proposition}
\begin{proof}
The first statement is Lemma~\ref{lem:pencil}(a). For the converse, we check that every line of $V$ is covered exactly once. A line of $V$ lies in $F$, or meets $F$ in exactly one point, or is disjoint from $F$. A line meeting $F$ in a point has a one-dimensional image in $V/F$, so it lies in exactly one $H_i$.

Lines of $F$ are covered exactly once by the $5$ common planes and by $\K_1,\K_2,\K_3$ (via $T_1,T_2,T_3$). A line of $H_i$ outside $F$ is a block line of $\K_i$, which is already covered, or free in $H_i$. A free line is covered by a transversal block iff it is one of the three lines through $P$ of an intercalate, and by Lemma~\ref{lem:pencil}(c),(d) an intercalate with parameters $(P,x,y)$ is the plane $\langle P,e_1+x,e_2+y\rangle$, whose three lines through $P$ are the three lines listed. There are $\sum_{P\in F}(6+2j(P))=6\cdot31+2\cdot35=256$ free lines of $H_i$ outside $F$, and $256$ intercalates, each contributing one line to each $H_i$. Hence if every free line is used exactly once, no block line of $\K_i$ is used. The lines disjoint from $F$ are the cells of the Cayley table, and they are covered exactly once iff the intercalates partition the table. Hence the condition says exactly that every line of $V$ is covered once.
\end{proof}

\begin{remark}\label{rem:pair}
For $\D|_{H_1}$ and $\D|_{H_2}$ it is necessary that they have the same $5$ blocks in $F$ and that $T_1\cap T_2=\emptyset$. Then $T_3$ is the complement and automatically has $5-j(P)$ lines through each $P$. From the point of view of a single hyperplane $H$ and a $5$-subspace $F\le H$ this gives a necessary condition for a hyperplane section to extend to a design: the $80$ free lines of $H$ inside $F$ must split into two sets of $40$ lines, with $5-j(P)$ lines of each through every point $P\in F$, each of which is the trace of a hyperplane section on one of the other two hyperplanes through $F$ with the same $5$ blocks in $F$. We have not been able to test this condition systematically.
\end{remark}

\subsection{Global pattern counts}
For a line spread $\mathcal S$ of $\PG(5,2)$ let $n_k(\mathcal S)$ be the number of solids containing exactly $k$ lines of $\mathcal S$.

\begin{lemma}[cf.\ {\cite[Sect.~3.C]{MT09}}]\label{lem:ident}
For every line spread $\mathcal S$ of $\PG(5,2)$, $\sum_kn_k=651$, $\sum_kk\,n_k=735$ and $\sum_k\binom k2n_k=210$. For a geometric spread, $(n_0,\dots,n_5)=(0,630,0,0,0,21)$.
\end{lemma}
\begin{proof}
Count solids, pairs (line of $\mathcal S$, solid through it) and triples (two lines of $\mathcal S$, solid containing both): $651=\gauss642$, $21\cdot\gauss422=735$, $\binom{21}2=210$, since two disjoint lines span a unique solid. In a geometric spread every solid contains $0$, $1$ or $5$ spread lines \cite[Sect.~2.3]{Kie22}. Then $10\,n_5=210$, $n_1=735-5\cdot21=630$ and $n_0=0$.
\end{proof}

Let $n_t$ be the number of $5$-subspaces of $V$ with pattern $t\in\{t_1,\dots,t_5\}$ (Proposition~\ref{prop:patterns}), and $N_j=\sum_Pn_j(\Der_P\D)$ (Lemma~\ref{lem:ident}).

\begin{proposition}\label{prop:glob}
$\sum_Fm_j(F)=N_j$ for all $j$. Therefore $n_{t_1}=N_5$ \textup{(}the number of $\beta$-flats\textup{)}, $n_{t_2}=N_4$, $2n_{t_3}+n_{t_4}=N_3$ and $\sum_tn_t=2667$. These are all linear relations, and the system has a solution in non-negative integers iff
\[
\textstyle\sum_Pw(\Der_P\D)\le2667,\qquad w(\mathcal S)=\tfrac12n_3(\mathcal S)+n_4(\mathcal S)+n_5(\mathcal S).
\]
This inequality is the sum over all $F$ of the inequalities $\tfrac12m_3(F)+m_4(F)+m_5(F)\le1$, whose left side takes the following values on the five patterns.
\begin{center}\small
\begin{tabular}{@{}l|ccccc@{}}
pattern & $t_1$ & $t_2$ & $t_3$ & $t_4$ & $t_5$\\\hline
$\tfrac12m_3+m_4+m_5$ & $1$ & $1$ & $1$ & $\tfrac12$ & $0$
\end{tabular}
\end{center}
\end{proposition}
\begin{proof}
A point $P\in F$ with $j_F(P)=j$ corresponds to a solid $F/P$ of $V/P$ containing $j$ lines of $\Der_P\D$, and every solid of $V/P$ arises from exactly one $F\ni P$, which gives $\sum_Fm_j(F)=N_j$. The relations for $j=0,1,2$ follow from those for $j=3,4,5$ and the identities of Lemma~\ref{lem:ident} and Proposition~\ref{prop:patterns}. The five pattern vectors span a $4$-dimensional space with the single relation $t_3+t_5=2t_4=(10,36,14,2,0,0)$. The solutions are $n_{t_1}=N_5$, $n_{t_2}=N_4$, $n_{t_4}=N_3-2n_{t_3}$, $n_{t_5}=2667-N_5-N_4-N_3+n_{t_3}$, and an integer $n_{t_3}$ with $\max(0,N_3+N_4+N_5-2667)\le n_{t_3}\le N_3/2$ exists iff $N_5+N_4+N_3/2\le2667$. Finally $\tfrac12m_3+m_4+m_5$ equals $1,1,1,\tfrac12,0$ on $t_1,\dots,t_5$.
\end{proof}

\begin{fact}\label{fact:w}
Every line spread $\mathcal S$ of $\PG(5,2)$ satisfies $w(\mathcal S)\le21$, as a computation over all $131\,044$ isomorphism types \cite{MT09} shows. Equality holds for exactly $8$ types, among them the Desarguesian spread ($n_5=21$, Lemma~\ref{lem:ident}). Thus $\sum_Pw(\Der_P\D)\le127\cdot21=2667$ for every assignment of spread types to the points, and the system of Proposition~\ref{prop:glob} is always solvable. Hence the linear relations of Proposition~\ref{prop:glob} between the numbers of $5$-subspace patterns and the types of the derived spreads cannot by themselves exclude a Fano design.
\end{fact}

\begin{proposition}\label{prop:lp}
Let $\mathcal P\subseteq\mathbb R^{11811}$ be the convex hull of the incidence vectors of all sets of planes of $\PG(6,2)$ pairwise meeting in at most a point, and $u=\frac1{31}\mathbf1$. The following are equivalent:
\begin{enumerate}[label=(\roman*)]
\item no Fano design exists,
\item $u\notin\mathcal P$,
\item $\sum_Bz_B\le380$ for all $z\in\mathcal P$.
\end{enumerate}
\end{proposition}
\begin{proof}
Every line lies in $31$ planes, so $u$ satisfies all line constraints with equality, and $\sum_Bu_B=11811/31=381$. (i)$\Rightarrow$(iii): a set $\mathcal C$ of planes pairwise meeting in at most a point covers $7|\mathcal C|\le2667$ lines, so $|\mathcal C|\le381$, with equality only for a design. (iii)$\Rightarrow$(ii) is clear. (ii)$\Rightarrow$(i), by contraposition: if $\D$ is a design, then $u\in\mathcal P$. Indeed, $\GL(7,2)$ is transitive on planes, so the average of the incidence vectors of $g\D$, $g\in\GL(7,2)$, is $u$.
\end{proof}

The polytope $\mathcal P$ is the integer hull of the set packing formulation of the instance \texttt{cdc7-4-3-2} of MIPLIB~2017 \cite{GHG+21} mentioned in the introduction.

\begin{corollary}
If a linear inequality $a\cdot z\le b$ is valid on $\mathcal P$ and violated by $u$, then $\sum_Bz_B\le380$ for all $z\in\mathcal P$. Hence a proof of non-existence by linear inequalities is a proof that every set of planes pairwise meeting in at most a point has at most $380$ elements.
\end{corollary}
\begin{proof}
Since $\mathcal P$ is invariant under $\GL(7,2)$, averaging over $\GL(7,2)$ gives the valid inequality $\bar a\sum_Bz_B\le b$ with $\bar a=\sum_Ba_B/11811$. As $0\in\mathcal P$, we have $b\ge0$, and $a\cdot u=381\bar a>b$ gives $\bar a>0$ and $\sum_Bz_B\le b/\bar a<381$ on $\mathcal P$. The vertices of $\mathcal P$ are incidence vectors of packings and have integral coordinate sums, so $\sum_Bz_B\le380$ on $\mathcal P$.
\end{proof}

\section{Conclusion and open problems}\label{sec:open}

For $q=2$ all bounds hold for abstract hyperplane sections, which exist (Appendix~\ref{app:section}), so a non-existence proof of the $q$-analog of the Fano plane has to use more than hyperplane sections and block counts in $5$-subspaces (Section~\ref{sec:dh}).

\begin{problem}
Improve the constant $41$ of Theorem~\ref{thm:41}. Every bound $|\X_H|\ge t$ transfers to designs via Lemma~\ref{lem:griesmer}, and for instance $t=56$ would give $|\N|\ge116$ and $|\A|\le11$.
\end{problem}

\begin{problem}
What is the maximum number of $\beta$-flats in a hyperplane section? By Theorem~\ref{thm:41} it is at most $22$, and $15$ is attained (Appendix~\ref{app:section}).
\end{problem}

Combining the hyperplane sections in the three hyperplanes through a $5$-subspace (Proposition~\ref{prop:pencil}, Remark~\ref{rem:pair}) may give better bounds on $|\X_H|$ or on the number of $\beta$-flats.

\appendix
\section{Programs and cross-checks}\label{app:repro}
All programs are written in C++. The hyperplane sections of Observation~\ref{obs:dh} were found by exact integer linear programming. We describe the main checks.

For Section~\ref{sec:tools}, the $15$ totally isotropic lines of $W(2)$, its spreads and their pairwise intersections were checked by enumeration. There are $6$ spreads, and any two of them share exactly one line (Lemma~\ref{lem:W2}).

For Section~\ref{sec:polar}, the facts on a nondegenerate alternating form on $\F_2^6$ were checked by direct enumeration. Every plane is TI ($135$ planes) or has rank $2$ with a unique radical point ($1\,260$ planes). For Lemma~\ref{lem:parity}, we enumerated the $56$ line spreads of $\PG(3,2)$ for each of the $28$ nondegenerate alternating forms on $\F_2^4$. A spread contains $1$, $3$ or $5$ totally isotropic lines, in $30$, $20$ and $6$ cases respectively.

For Fact~\ref{fact:local}, we enumerated the $18\,958\,856$ sets of five pairwise disjoint lines of $\PG(4,2)$, which form $10$ orbits under $\GL(5,2)$, and checked the extension property for every orbit using the line spreads of $\PG(5,2)$ of \cite{MT09}.

For Fact~\ref{fact:w}, we computed $n_3(\mathcal S)$, $n_4(\mathcal S)$ and $n_5(\mathcal S)$ for each of the $131\,044$ isomorphism types of line spreads of $\PG(5,2)$ of \cite{MT09} and evaluated $w(\mathcal S)$. The maximum is $21$, attained for exactly $8$ types.

\section{An explicit hyperplane section}\label{app:section}
Identify a point of $H=\F_2^6$ with the integer $v=\sum_{i=0}^5v_i2^i\in\{1,\dots,63\}$ formed by its coordinates, so that addition is the bitwise XOR. A triple $(a,b,c)$ of linearly independent vectors stands for the plane $\langle a,b,c\rangle$ with the seven points $a,b,c,a\oplus b,a\oplus c,b\oplus c,a\oplus b\oplus c$. The following $45$ planes form a hyperplane section $\K$. The planes pairwise meet in at most a point, and every point lies on exactly $5$ of them. A second program, written independently of the search code, confirms this.
\begin{center}\small
\setlength{\tabcolsep}{6pt}
\begin{tabular}{@{}ccccc@{}}
$(1,2,12)$ & $(3,8,37)$ & $(5,19,43)$ & $(7,24,34)$ & $(13,17,37)$\\
$(1,4,24)$ & $(3,9,33)$ & $(5,26,42)$ & $(7,25,40)$ & $(14,20,37)$\\
$(1,6,18)$ & $(3,16,41)$ & $(5,27,41)$ & $(8,20,36)$ & $(14,21,38)$\\
$(1,8,22)$ & $(3,17,32)$ & $(6,10,32)$ & $(8,21,39)$ & $(14,22,36)$\\
$(1,10,16)$ & $(4,18,34)$ & $(6,11,51)$ & $(9,20,35)$ & $(14,23,39)$\\
$(2,8,49)$ & $(4,19,40)$ & $(6,24,33)$ & $(9,21,32)$ & $(15,20,34)$\\
$(2,9,52)$ & $(4,26,32)$ & $(6,25,43)$ & $(12,16,33)$ & $(15,21,33)$\\
$(2,16,37)$ & $(4,27,35)$ & $(7,10,35)$ & $(12,17,34)$ & $(15,22,35)$\\
$(2,17,45)$ & $(5,18,33)$ & $(7,11,48)$ & $(13,16,38)$ & $(15,23,32)$\\
\end{tabular}
\end{center}
The section has exactly $15$ $\beta$-flats, so $|\X_H|=48$. They are the kernels of the linear forms $x\mapsto\mathrm{parity}(x\,\&\,f)$ for
\begin{align*}
f&=2,8,10,20,22,28,30,32,34,40,42,52,54,60,62,\\
\text{with focal points }&28,18,26,9,21,15,6,1,29,19,27,8,20,14,7
\end{align*}
in this order. These $15$ focal points are exactly the points of a solid, so $\A'$ spans only a $4$-dimensional space and the bound of Lemma~\ref{lem:nonspan}(b) is attained. The three $5$-subspaces through this solid are $\beta$-flats (with $f=22,32,54$), and each of them meets $\X$ in $16$ points.

\section*{Declaration of generative AI use}
The author used Claude Sonnet 5.5 (Anthropic) to improve the language and wording of the manuscript and to check proofs for errors, inconsistencies, and notation clashes. All mathematical statements and proofs were obtained and verified by the author, who takes full responsibility for the content of this paper.

\end{document}